\documentclass[11pt,a4paper]{amsart}
\begin{document}

\newcommand{\ea}{\mbox{{\bf a}}}
\newcommand{\eu}{\mbox{{\bf u}}}
\newcommand{\ueu}{\underline{\eu}}
\newcommand{\ueo}{\overline{u}}
\newcommand{\oeu}{\overline{\eu}}
\newcommand{\ew}{\mbox{{\bf w}}}
\newcommand{\ef}{\mbox{{\bf f}}}
\newcommand{\eF}{\mbox{{\bf F}}}
\newcommand{\eC}{\mbox{{\bf C}}}
\newcommand{\en}{\mbox{{\bf n}}}
\newcommand{\eT}{\mbox{{\bf T}}}
\newcommand{\eL}{\mbox{{\bf L}}}
\newcommand{\eR}{\mbox{{\bf R}}}
\newcommand{\eV}{\mbox{{\bf V}}}
\newcommand{\eU}{\mbox{{\bf U}}}
\newcommand{\ev}{\mbox{{\bf v}}}
\newcommand{\eve}{\mbox{{\bf e}}}
\newcommand{\uev}{\underline{\ev}}
\newcommand{\eY}{\mbox{{\bf Y}}}
\newcommand{\eK}{\mbox{{\bf K}}}
\newcommand{\eP}{\mbox{{\bf P}}}
\newcommand{\eS}{\mbox{{\bf S}}}
\newcommand{\eJ}{\mbox{{\bf J}}}
\newcommand{\eB}{\mbox{{\bf B}}}
\newcommand{\leb}{\mathcal{ L}^{n}}
\newcommand{\eI}{\mathcal{ I}}
\newcommand{\eE}{\mathcal{ E}}
\newcommand{\hen}{\mathcal{ H}^{n-1}}
\newcommand{\eBV}{\mbox{{\bf BV}}}
\newcommand{\eA}{\mbox{{\bf A}}}
\newcommand{\eSBV}{\mbox{{\bf SBV}}}
\newcommand{\eBD}{\mbox{{\bf BD}}}
\newcommand{\eSBD}{\mbox{{\bf SBD}}}
\newcommand{\ecs}{\mbox{{\bf X}}}
\newcommand{\eg}{\mbox{{\bf g}}}
\newcommand{\paromega}{\partial \Omega}
\newcommand{\gau}{\Gamma_{u}}
\newcommand{\gaf}{\Gamma_{f}}
\newcommand{\sig}{{\bf \sigma}}
\newcommand{\gac}{\Gamma_{\mbox{{\bf c}}}}
\newcommand{\deu}{\dot{\eu}}
\newcommand{\dueu}{\underline{\deu}}
\newcommand{\dev}{\dot{\ev}}
\newcommand{\duev}{\underline{\dev}}
\newcommand{\weak}{\rightharpoonup}
\newcommand{\weakdown}{\rightharpoondown}
\newtheorem{rema}{Remark}[section]
\newtheorem{thm}{Theorem}[section]
\newtheorem{lema}{Lemma}[section]
\newtheorem{prop}{Proposition}[section]
\newtheorem{defi}{Definition}[section]
\newtheorem{exempl}{Example}[section]
\newtheorem{opp}{Open Problem}[]
\renewcommand{\contentsname}{ }
\newenvironment{rk}{\begin{rema}  \em}{\end{rema}}
\newenvironment{exemplu}{\begin{exempl}  \em}{\end{exempl}}

\title{QUASICONVEXITY VERSUS GROUP INVARIANCE}
\author{MARIUS BULIGA \email{marius.buliga@imar.ro}}
\address{Institute of Mathematics of the Romanian Academy \\ 
 PO BOX 1-764, RO 70700, Bucharest}

\begin{abstract}
The lower invariance  
under a given  arbitrary group of diffeomorphisms extends the notion 
of quasiconvexity. 
The non-commutativity 
of the group operation (the function composition) modifies the 
classical equivalence between lower semicontinuity and 
quasiconvexity.

In this context  null lagrangians are particular
cases of integral invariants of the group. 

\vspace{1.cm}

\end{abstract}

\maketitle

\markboth{Marius Buliga}{Quasiconvexity Versus Group Invariance}

%\vspace{2cm}

%\newpage

\section{Introduction}

\subsection{First notations}

In this paper\footnote{Lecture held on Feb. 22 at the Mathematical Institute, 
Oxford, Applied Analysis and Mechanics Seminars,Hilary 
Term 1999.} $\Omega \subset {\mathbb R}^{n}$ is an open
 bounded set, with smooth boundary and  
$id$ is the identity map of ${\mathbb R}^{n}$ $id (x) = x$. 
$A \subset \subset B$ means that  $A$ is compactly included in  
$B$.

 $Diff^{\infty}_{0} = Diff^{\infty}_{0}(R^{n})$ is the group of 
all $C^{\infty}$ diffeomorphisms with compact support in  
${\mathbb R}^{n}$:
$$Diff^{\infty}_{0}(R^{n}) \ = \ \left\{ \phi \in C^{\infty}
({\mathbb R}^{n},
{\mathbb R}^{n}) \mbox{ : } \exists \ 
 \phi^{-1} \in C^{\infty}({\mathbb R}^{n},{\mathbb R}^{n}) , \right.
$$ 
$$ \left. supp \ (\phi - id) \subset \subset {\mathbb R}^{n} 
\right\} . $$

For any open set $\Omega \subset \subset R^{n}$ and any 
subgroup $G \subset Diff^{\infty}_{0}$ we define 
$$G(\Omega) \ = \ \left\{ \phi \in  Diff^{\infty}_{0} \mbox{ : } 
supp  \ (\phi - id) \ \subset  \Omega \right\} \ \ . $$
Notice that $G(\Omega)$ is a group under the operation "." of 
functions composition.

The main goal of this paper is to study the 
lower semicontinuity  of an  integral functional having the form:
$$I(\eu, \Omega) \ = \ \int_{\Omega} W(x, \eu(x), \nabla \eu (x)) 
\mbox{ d}x$$ 
under inner (or left) and outer (or right)
variations in a group of diffeomorphisms $G(\Omega)$. 
An inner (or left) variation of $I$ consists in the replacement of the 
argument $\eu$ by $\eu . \phi$ with $\phi \in G(\Omega)$. Similarly, 
an outer (or right) variation consists in the replacement of $\eu$ by 
$\phi . \eu$. Off course, the argument $\eu$ belongs to a 
space of mappings $X$ which is stable under inner or outer variations. 

Lower semicontinuity results are important tools in the study of 
existence and regularity of critical points for (integral) functionals.

\subsection{Some problems involving groups of diffeomorphisms}

In this section we shall give three examples of problems involving 
groups of diffeomorphisms and critical points of integral functionals.

\subsubsection{Least action principles in nonlinear mechanics}

The standard space of configuration of
a fluid lying in a vessel $\Omega$,
under adherence conditions on the
vessel's wall $\paromega$, is $Diff_{0}^{\infty}(\Omega)$.
Moreover, if the fluid is incompressible, then the space of
states is $ Diff_{0}^{\infty}(dx)(\Omega)$, namely the group 
of volume-preserving diffeomorphisms with support in $\Omega$.

Arnold  \cite{ar1} first showed that hydrodynamics of an 
ideal fluid can be formulated in the frame of volume preserving 
diffeomorphisms: the evolution of an ideal fluid is a geodesic 
(i.e. minimizer of an integral action functional) in the 
group of volume preserving diffeomorphisms. This work has been 
developed in papers like Arnold \& Khesin \cite{ak1}, Ebin \& Marsden 
\cite{em} or Shnirelman \cite{sn}.

Of a related nature is the problem of evolution of an elastic
body. Generally, we seek for critical points of an action 
functional $$A(\eu_{t}) \ = \ 
\int_{\Omega \times [0,T]} L(x,t, \eu(x,t), \nabla \eu(x,t)) 
\mbox{ d}x\mbox{ d}t$$ under admissible variations of 
the form 
$$(t, \eu( \cdot, t)) \in [0,T] \times X(\eu^{0}(t)) 
\mapsto (t, \ev(\cdot,t)) \in [0,T] \times X(\eu^{0}(t)) \ \
$$
such that $\ev(\cdot, 0) = \eu (\cdot, 0)$, 
$\ev(\cdot, T) = \eu (\cdot, T)$, $\dot{\ev}(\cdot, 0) =
\dot{\eu} (\cdot, 0)$, $\dot{\ev}(\cdot, T) = \dot{\eu} 
(\cdot, T)$. This problem can be reformulated by considering 
variations of the form 
$$(t, \eu( \cdot, t)) \in [0,T] \times X(\eu^{0}(t)) 
\mapsto (t, \eu . \phi (\cdot, t)) \in [0,T] \times X(\eu^{0}(t)) \ \
$$
for all $$ \phi \in Diff^{\infty}_{0}(\Omega \times [0,T]) \ , \ 
\phi(x,t) \ = \ (\phi_{t}(x), t) \ , \ \phi_{t} \in 
Diff^{\infty}_{0}(\Omega) \ \ . $$
This set of diffeomorphisms of $\Omega \times [0,T]$ is a
subgroup of $Diff^{\infty}_{0}(\Omega \times [0,T])$.

\subsubsection{Critical points turned into local minima}

Another interesting problem is to find all critical points of 
an integral functional by variational methods. 
The direct method gives access only 
to global minimizers. Zhang \cite{z1}, \cite{z2},  Sivaloganathan 
\cite{s1} show that  for  strictly quasiconvex 
potential the critical points of the associated functional 
are local (the word has various precise meanings) minimizers. 

The connection of the critical point problem with
diffeomorphisms groups is made by the following  proposition:

\begin{prop}
Suppose that $W:  GL_{n}(R)
\rightarrow  R$ is a $C^{2}$ potential and   
$\eu \in C^{3}(\Omega, R^{n})$ is a local minimizer of the 
functional 
$$I(\ev) \ = \ \int_{\Omega} W(\nabla \ev) \mbox{ d}x $$
in the class 
$$\eu Diff_{0}^{\infty}(\Omega) \ = \ \left\{ \eu . \phi 
\mbox{ : }  \phi \in Diff_{0}^{\infty}(\Omega) \right\} \ \ . $$ 
Then $\eu$  is a critical point of $I$. 
\label{pcr}
\end{prop}

\begin{proof}
Let us consider, for a given but arbitrary 
$\eta \in C^{\infty}_{c}(\Omega, R^{n})$, 
the one parameter flow $t \mapsto \phi_{t}$ generated by $\eta$. 
We have then $I(\eu) \leq I(\eu . \phi_{t}^{-1})$ for small 
$\mid t \mid$. Therefore:
$$\frac{\partial I}{\partial t}(\eu . 
\phi_{t}^{-1})_{|_{t = 0}} \ = \ 0 \ \ .$$
The latter equality means that for any $\eta \in 
C^{\infty}_{0}(\Omega, R^{n})$ we have: 
$$\int_{\Omega} \left\{ W(\nabla \eu) \ div \ \eta \ - 
\ \frac{\partial W}{\partial \eF_{ik}} 
(\nabla \eu) \eu_{i,l}\eta_{l,k} \right\} \mbox{ d}x \ = 
\ 0 \ \ . $$
This is equivalent to: 
$$ \int_{\Omega} \left\{ W(\nabla \eu) \eta_{i} - 
\frac{\partial W}{\partial \eF_{lk}} 
(\nabla \eu) \eu_{l,i}\eta_{i}\right\}_{,i} \ - $$ 
$$ \left\{ W(\nabla \eu)_{,i}\eta_{i} + 
\left(\frac{\partial W}{\partial \eF_{lk}} (\nabla \eu) 
\eu_{l,i}\right)_{,l}\eta_{i} \right\} 
\mbox{ d}x \ = \ 0 \ . $$ 
By the Divergence theorem we have: 
$$\int_{\Omega} \left\{ \left(\frac{\partial W}{\partial 
\eF_{lk}} (\nabla \eu)\right)_{,l} 
\eu_{l,i} \eta_{i} \right\} \mbox{ d}x \ = \ 0 \ . $$ 
We conclude that 
$$\left(\nabla \eu \right)^{T} \ 
div \left(\frac{ \partial W}{\partial \eF_{lk}} 
(\nabla \eu) \right) \ = \ 0 \ .$$
Because $\nabla \eu(x)  \in GL_{n}(R)$ it follows that 
$$div \left(\frac{ \partial W}{\partial \eF_{lk}} 
(\nabla \eu) \right) \ = \ 0 \ .$$
\end{proof}

We see from proposition \ref{pcr} that the class of all
(sufficiently regular) $\eu$ 
with the property: 
\begin{equation}
I(\eu) \leq I(\eu . \phi) \mbox{ for all } \phi 
\in Diff_{0}^{\infty}(\Omega) 
\label{pgob}
\end{equation} 
is included in the class of critical points of $I$.  
From the proof of the same proposition we notice that any critical 
point $\eu$ of $I$ is critical in the class 
$\eu . Diff_{0}^{\infty}(\Omega)$ in the sense that for any $\eta 
\in C^{\infty}_{c}(\Omega, R^{n})$ which 
generates the one parameter flow $\phi_{t}$  we have: 
\begin{equation}
\frac{d}{dt} I(\eu . \phi_{t}^{-1} )_{|_{t = 0}} \ =  \ 0 
\label{cgr}
\end{equation}  
A natural question is: does any critical point $\eu$ have the property 
\eqref{pgob}?

Consider, as an example,  the functional 
$$I(\theta) \ = \ \int_{0}^{L} \frac{1}{2} k \mid \theta' \mid^{2} + 
\lambda \eF (\theta) \mbox{ d} s \ \ ,$$
where $\eF$ is an arbitrary $C^{1}$ function. 

The group $Diff_{0}^{\infty} (0,L)$ is nothing but: 
$$Diff_{0}^{\infty} (0,L) \ = \ \left\{ \phi \in C^{\infty}((0,L),
(0,L)) \mbox{ : } \phi' > 0 \ , \  supp \ (\phi - \  id) \subset 
\subset (0,L) \right\} \ . $$
We introduce now the functional:
$$I(\phi; \theta) \ = \ I(\theta . \phi^{-1}) \ = \ 
\int_{0}^{L} \frac{1}{2} k \mid \theta' \mid^{2} \frac{1}{\phi'} + 
\lambda \eF (\theta) \phi' \mbox{ d}s \ \ . $$  
It is easy to see that $I( \cdot ; \theta)$ is convex on  
$Diff_{0}^{\infty} (0,L)$. 

Consider now $\theta$ a critical point of $I$. 
The function $id$ is then a critical point of $I(\cdot; \theta)$
 in $Diff_{0}^{\infty} (0,L)$. 
We deduce from the convexity of $I(\cdot; \theta)$
that $id$ is a global minimum, therefore  it satisfies the relation 
\eqref{pgob}.

\subsubsection{The invariance problem}

Let $X$ be a class of functions from $R^{n}$ to $R^{m}$ 
and $G$ a topological semigroup of diffeomorphisms of 
$R^{n}$ such that $X.G \ = \ X$.

\begin{defi}
The functional $I: X \rightarrow R$ is $G$ left  
invariant if for any $\eu \in X$ and $\phi \in G$ we 
have the equality $I(\eu . \phi) \ = \ I(\eu)$. 
\label{dinva}
\end{defi}

 Such invariants play a central role in continuum mechanics. 
Indeed, let us consider $X \ = \ C^{3}(\Omega, R^{k})$, 
$G = Diff_{0}^{\infty}(\Omega)$ and 
$$I(\eu; \Omega) \ = \ \int_{\Omega} L(\nabla \eu (x))
\mbox{ d}x \ \ . $$  
We prove in this paper the following 

\begin{prop}
Let $L: M^{n, n} \rightarrow R$ be a $C^{2}$ map. We have
equivalence between the statements:
\begin{itemize}
\item[(i)] $I$ is a  $Diff_{0}^{\infty}(\Omega)$ invariant,
\item[(ii)] $L$ is a null lagrangian. 
\end{itemize}
\label{pinvlag}
\end{prop}

\section{Outline}

We address in this paper the problem of proving lower
semicontinuity of integral functionals defined over groups 
of diffeomorphisms. 

After a section of preliminaries we introduce the notion of 
left and right lower  invariance of a 
mapping with respect to a group 
and prove that for mappings defined on $GL_{n}(R)$, 
quasiconvexity in the sense of Morrey is equivalent to the 
left lower invariance under the group of diffeomorphisms 
with compact support. 

The lower invariance  
under a given  arbitrary group $G$ of diffeomorphisms extends 
the notion of quasiconvexity.  
The non-commutativity 
of the group operation (the function composition) modifies the 
classical equivalence between lower semicontinuity and 
quasiconvexity. We introduce a new notion of semicontinuity, 
named $G$ left (or right) lower semicontinuity.

 The main results of this paper (theorems 
\ref{t2} and  \ref{t3}) show that  if 
 the integral functional $I$ is $G$ left lower semicontinuous 
 then  the potential  $W$ is $G$ left lower invariant; also   
 if  $W$ is $G$  right lower invariant    then  $I$ is $G$ 
right lower semicontinuous. 

In this context  null lagrangians are particular
cases of integral invariants of the group. We generally find that the 
only homogeneous integral functionals which are  weak * 
continuous under inner variations in a group  of 
diffeomorphisms are the constant ones.

\section{Preliminaries}

\subsection{Notations}

The set $W^{1,\infty}(\Omega)$ is the 
 Sobolev space of  $L^{\infty}(\Omega, {\mathbb R}^{n})$ 
functions with the first derivative essentially bounded.  
$\| \cdot \|_{1,\infty}$ is the usual norm in   
$W^{1,\infty}(\Omega, {\mathbb R}^{n})$.   
For any  $\eu \in W^{1,\infty}(\Omega, {\mathbb R}^{n})$ 
the gradient $\nabla \eu$ is identified with the approximate 
gradient of $\eu$. Therefore the equalities involving gradients  
are almost everywhere (abbreviated "a.e.").

$\mathcal{ A}$ is the group of affine homothety-translations. Any 
element of $\mathcal{ A}$ has the form:  
$$ \alpha(x_{0},y_{0},\epsilon)(x) \ = \ \ef (x) \ = \ x_{1} + 
\epsilon (x -x_{0}) \ \ , \ \ x_{0}, 
x_{1}  \  \in \ {\mathbb R}^{n} \ , \ \epsilon \ > \ 0 \ \ .$$ 
We consider on $\mathcal{ A}$ the punctual convergence of functions 
defined on  ${\mathbb R}^{n}$ with values in ${\mathbb R}^{n}$.  
$GL_{n}(R) \subset {\mathbb R}^{n \times n}$ is the
multiplicative group of  all 
invertible, orientation preserving,
 matrices, i.e the set of all $\eF$ such that  
$det \ \eF \ > 0$.

We shall use in the paper the affine space 
$$W^{1,\infty}_{id}({\mathbb R}^{n}) \ = \ \left\{ 
\phi  \in W^{1,\infty}_{loc}({\mathbb R}^{n}, {\mathbb R}^{n}) 
\mbox{ : } supp \  (\phi \ - \ id ) \subset \subset 
{\mathbb R}^{n} \right\} \ \ .$$

\subsection{Basic definitions and properties}

 $G$ is a set  of functions from  
${\mathbb R}^{n}$ to  ${\mathbb R}^{n}$, 
which satisfies the following axioms:
\begin{enumerate}
\item[A1/] $(G,.)$ is a semigroup with  
the function composition operation ".";
\item[A2/] $G  \ \subset \ W^{1,\infty}_{id}({\mathbb R}^{n}) 
\cap C^{1}({\mathbb R}^{n},{\mathbb R}^{n})$ ;
\item[A3/] the following action is well defined:
$$A: \mathcal{ A} \times G \rightarrow G \ \ , \ \ A( \ef, \phi )
 \ = \ \ \ef . \phi . \ef^{-1} \ \ . $$ 
\end{enumerate}
Further on we shall suppose that $G$ acts transitively on 
${\mathbb R}^{n}$. We did not included this
statement among the axioms because all the results from the 
paper hold without the transitivity assumption, but in a more 
involved form. The same remark is true if we suppose only that 
for any $x \ \in \Omega$ the orbit $$G(\Omega)(x) \ = \ 
\left\{ \phi(x) \mbox{ : } \phi \in G(\Omega) \right\}$$
 is dense in $\Omega$. For any open set $E \subset \subset
R^{n}$ the set $G(E)$ is defined further (definition
\ref{dunu}).

\begin{defi} For any open set $E \subset {\mathbb R}^{n}$ 
we define 
$$G(E) \ = \ \left\{ \phi \in G \mbox{ : } supp \ (\phi - id) 
\subset \subset E \right\} \ \ . $$

For any $x_{0} \in {\mathbb R}^{n}$  
the first order jet of $G$ in $x_{0}$ is: 
$$J^{1}(x_{0},G) \ = \ \left\{  
\nabla \phi (x_{0}) \mbox{ : } \phi \in G \right\} \ \ . $$
\label{dunu}
\end{defi}

\begin{defi}
Let $(\phi_{h})_{h} \in G(\Omega)$ be a sequence and 
$\phi$ an element of  $W^{1,\infty}_{id}(\Omega, 
{\mathbb R}^{n})$. We say that $\phi_{h}$  converges 
to $\phi$  
if the sequence $\phi_{h}$ converges $W^{1,\infty}$ 
weak* to $\phi$.

 $G^{1,\infty}(\Omega)$  is the closure of $G$ in 
$W^{1,\infty}_{id}(\Omega, {\mathbb R}^{n})$ with 
respect the strong convergence. That is $G^{1,\infty}
(\Omega)$ is the space of all $\eu$ which 
can be obtained as limit points of strong convergent sequences
$\phi_{h}$, $\phi_{h} \in G(\Omega)$.   
\label{def2}
\end{defi} 

\begin{rk}
The action $A$ is continuous. The  operation "." 
is continuous in each argument.
\label{r0}
\end{rk}

In the following lemmas we collect some elementary 
facts connected to the convergence  or to the algebraic 
structure previously introduced.

\begin{lema}
\begin{enumerate}
\item[1/] Let $\Omega \subset B(0,R)$. Then 
$G(\Omega) \subset G(B(0,R))$. 
\item[2/] $G^{1,\infty}(\Omega) \subset C^{0,1}(B(0,R))$.
 For any sequence 
$\phi_{h} \in G^{1,\infty}(\Omega)$ such that $\phi_{h} 
\stackrel{b}{\rightarrow} \phi$ there exists a
subsequence which converges uniformly to $\phi$. 
\item[3/] $G^{1,\infty}(\Omega)$ is a semigroup. The 
composition operation "." is continuous in each
argument.   
\end{enumerate}
\label{ltool}
\end{lema}

\begin{lema}
Let $A, B$ be non empty open subsets of ${\mathbb R}^{n}$.
 If $A$ is bounded then there exists  
$\ef \in \mathcal{ A}$ such that the application $A(\ef, \cdot) :
 G(A) \rightarrow G(B)$ is well defined, 
injective and continuous. 
\label{p1}
\end{lema}

\begin{lema}
If $A$ and $B$ are two open disjoint sets then for any  
$\phi \in G(A)$, $\psi \in 
G(B)$ we have $\phi . \psi \ = \ \psi . \phi \ \in \ 
G(A \cup B)$ . 
\label{p2}
\end{lema}

The following proposition shows that the first order jet 
associated to the set $G$ and an arbitrary point $x \in {\mathbb R}^{n}$ 
 is a (semi) group which does not depend  on $x$.

\begin{prop}
There exists $J(G)$ sub semigroup of the multiplicative
 group $GL_{n}(R)$ 
such that for any 
$x_{0} \in {\mathbb R}^{n}$ we have $J^{1}(x_{0},G) \ = \ 
J(G)$.  If  $G$ is a group then  $J(G)$ is 
a group. 
\label{p3}
\end{prop}

\begin{proof} We first prove that $J^{1}(x_{0},G)$ is 
semigroup. We have $id \in G$, 
hence $I$, the identity matrix, belongs to $J^{1}(x_{0},G)$. 
Let us consider 
$R, S \in J^{1}(x_{0},G)$ and  
 $\phi, \psi \in G$ such that $R = 
\nabla \phi (x_{0})$, 
$S  =  \nabla \psi(x_{0})$. We define the translation 
$\ef \in \mathcal{ A}$: 
$\ef(x)  =  x + \psi(x_{0}) - x_{0}$. From A3/ we have 
$ \tilde{\phi} \ = \ \ef . \phi . \ef^{-1} \in G$, hence  from  
$\nabla \tilde{\phi} (\psi(x_{0}))  =  \nabla \phi (x_{0})$ 
and A2/  we infer that  
$$RS \ = \ \nabla \phi (x_{0}) \nabla \psi(x_{0}) \ = \ 
\nabla \tilde{\phi} (\psi(x_{0})) 
 \nabla \psi(x_{0}) \ = $$ $$ = \ \nabla ( \tilde{\phi} . 
\psi) (x_{0})  \ \in J^{1}(x_{0}, G) \ \ . $$
A simple argument based on A3/ shows that $J^{1}(x_{0},G)$ 
does not depend on  
 $x_{0} \in {\mathbb R}^{n}$. For a fixed, arbitrarily 
chosen,  $x_{0}$ we define $J(G) \ = \ 
J^{1}(x_{0},G)$. 

The proof of the fact that if $G$ is a group and $\eF 
\in J(G)$ then  
$\eF^{-1}$ exists and $\eF^{-1} \in J(G)$ is similar. 
\end{proof}

\begin{defi}
$W^{1,\infty}(G, \Omega)$ is the class of all $\eu \in 
W^{1, \infty}(\Omega, {\mathbb R}^{n})$ such that we have 
$\nabla \eu(x) \in J(G)$ a.e. in $\Omega$. 
\label{dclas}
\end{defi}

We describe now  two easy procedures of 
construction of groups satisfying the axioms 
A1/, A2/, A3/.

\begin{defi}
For any subgroup $M$ of  $GL_{n}(R)$ 
we define  the local group generated by $M$: 
$$[M] \ = \ \left\{ \phi \in Diff^{\infty}_{0} \mbox{ : } 
\forall \ 
  x \in {\mathbb R}^{n} \ \ 
\nabla \phi( x)  \  \in \  M  \  \right\} \ \ . $$ 
\label{dlocal}
\end{defi}
It is obvious that $[M]$ satisfies the axioms and 
$J([M]) \ = \ M$.  
We notice that the groups $G$ constructed in this way are 
determined by $J(G)$, that is~:
if  $J(G_{1})  = 
J(G_{2})$ then $G_{1}  =  G_{2}$.  This property justifies 
the name "local group".

\begin{defi}
For any semigroup (group) $G$ which satisfies 
the axioms  the completion of $G$ is defined by:
$$G^{c} \ = \ \left\{ \eF.\phi. \eF^{-1} \mbox{ : } \eF, \  
\eF^{-1} \in J(G) \ , \ \phi \in G 
\right\} \ \ . $$ 
\label{dcompl}
\end{defi}
Generally $G^{c}$ is  larger than $G$, but not always a  semigroup (group).

\begin{exemplu}
We obviously have $[GL_{n}(R)] \ = \ Diff^{\infty}_{0}$, 
therefore $J(Diff^{\infty}_{0}) =  GL_{n}(R)$.  
We have also $Diff^{\infty, c}_{0} \ = \ Diff^{\infty}_{0}$. 
\label{ex1}
\end{exemplu}
\begin{exemplu} Let us consider $Diff^{\infty}_{0}(dx)$, the 
subgroup of $Diff^{\infty}_{0}$ containing all 
volume preserving smooth diffeomorphisms with compact support. 
We have $Diff^{\infty}_{0}(dx)  =  [SL_{n}(R)] = 
Diff^{\infty, c}_{0}(dx)$. 
\label{ex2}
\end{exemplu}

\begin{exemplu} 

For any  $\eu: {\mathbb R}^{2n} \rightarrow {\mathbb R}^{2n}$ 
and $\omega$,the canonical symplectic 2-form on 
${\mathbb R}^{2n}$,  we 
denote by 
 $\eu^{*}(\omega)$  the transport of $\omega$. Let us define 
$Diff^{\infty}_{0}(\omega)$: 
$$Diff^{\infty}_{0}(\omega) \ = \ \left\{ \phi \in 
Diff^{\infty}_{0} \mbox{ : } 
\phi^{*}(\omega) \ = \ \omega \ \right\} \ \ .$$
The axioms A1/, A2/ A3/ are satisfied. We  have the equalities: 
$$J(Diff^{\infty}_{0}(\omega) ) \ = \ Sp_{n}(R) \ = \  
\left\{ \eF \in {\mathbb R}^{2n \times 2n} \mbox{ : }  
\eF \omega \eF^{T} \ =  \ \omega \ \right\} $$
and $Diff^{\infty, c}_{0}(\omega) \ 
= \ Diff^{\infty}_{0}(\omega) \ = \ [Sp_{n}(R)]$. 
\label{ex3}
\end{exemplu}

\begin{exemplu} 
Let us take a group $G$ which satisfies the axioms. The space 
of smooth loops $t \in S^{1} \mapsto \phi_{t} \in G$ 
can be embedded in the following group: 
$$LG \ = \ \left\{ \phi \in Diff^{\infty}_{0}({\mathbb R}^{n+1})
 \mbox{ : } \phi(t,x) \ = \ (t,\phi_{t}(x)) \ , \ \phi_{t} \in 
G \right\} \ \ . $$
Notice that $LG$ does not satisfy the transitivity assumption.
However, the results from this paper are true in this case, but 
with minor modifications which are left to the interested
reader. 
\label{exloop}
\end{exemplu}

\begin{exemplu}
Consider the class $\mathcal{ H}$ of hamiltonian diffeomorphisms
with compact support 
of ${\mathbb R}^{2n}$ (see Hofer, Zehnder \cite{10}). This is a 
group which satisfies the axioms, but it is not local. However, it is complete. 
\label{hamloop}
\end{exemplu}
For all the transitivity results needed in these example we 
refer to Michor \& Vizman \cite{8}.

%
%
%
%           Rezultate de semicontinuitate
%
%
%creat pe 29.01.98, destinat includerii in full_lecture.tex
%
%
%
%
\section{Lower invariance and semicontinuity}

The lower semicontinuity of functionals $I(\cdot; \Omega)$ 
defined over Sobolev spaces was systematically studied. 
Morrey \cite{9} introduced the notion of quasiconvexity and 
proved that $W^{1,\infty}$ weak * lower semicontinuity of 
$I(\cdot; \Omega)$ is equivalent to the quasiconvexity of 
the integrand $W$ in it's third variable, provided that 
$W$ is continuous. Acerbi \& Fusco \cite{2}, Ball \& Murat 
\cite{1"} improved this result and introduced the notion of 
$W^{1,p}$ quasiconvexity. Ball \cite{1}, \cite{1'}, considered 
a new condition, called polyconvexity, which implies
quasiconvexity, with important applications in nonlinear 
elasticity.

\subsection{Lower invariance and quasiconvexity}

There are several slightly different 
definitions of quasiconvexity. We prefer the one from Ball \cite{1'}: 
\begin{defi}
$W$ is quasiconvex in $(x_{0}, y_{0}, \eF) \in {\mathbb R}^{n}
 \times {\mathbb R}^{n} \times GL_{n}(R)$ if for 
any open bounded set $E \subset {\mathbb R}^{n}$ and any $\eta 
\in C^{\infty}( E, {\mathbb R}^{n})$ such 
that: 
\begin{enumerate}
\item[i)] $supp \ \eta \ \subset \subset E$ , 
\item[ii)] for any  $x \in E$ we have $\eF  + \nabla \eta (x)
 \in GL_{n}(R)$ (i.e. $det \ ( 
\eF  + \nabla \eta (x) )  \ > \ 0$) , 
\end{enumerate}
we have the inequality: 
\begin{equation}
\int_{E}W(x_{0}, y_{0}, \eF + \nabla \eta (y)) \mbox{ d}
 y \ \geq \ \mid E \mid 
W(x_{0}, y_{0}, \eF) \ \ . 
\label{qconv}
\end{equation}
\label{defqc}
\end{defi}

We could work in the followings 
  with a topological space $X$ of continuous 
functions from $R^{n}$ to $R^{m}$, such that 
$X . G \ = \ X$ and for any homothety-
translation $h:R^{n} \rightarrow R^{n}$ and any $\eu \in 
X$ we have $\eu . h \in X$.
Our model space will be 
$W^{1,\infty}(G, \Omega)$ (see definition \ref{dclas}); all 
results from this paper can be reformulated for a wide variety 
of spaces $X$ in an obvious way. We leave this for further 
applications.

The non-commutativity of the function composition forces us to consider 
a "left" and "right" variant of any further definition.
For the notions involving the word "right" we shall suppose 
that $G . X \ = \ X$ and  $X \subset 
W^{1,\infty}_{loc}(R^{m},R^{n})$. 

We introduce the following definition of lower (upper
respectively) invariance. In the next definition 
 we shall denote by $J(X)$ the first order jet of $X$ (supposing that 
it does not depend on $x$). 
From the condition   $X . G = X$ we derive  that  $J(X).J(G) \ = \ J(X)$. 

\begin{defi} Let us consider $x_{0}, y_{0} \in 
{\mathbb R}^{n}$ and $\eF \in J(X)$.

The function $W$ or the functional $I$   are  $G$ left lower 
invariant in $(x_{0}, y_{0}, \eF)$, and we shall write 
"$G$ L.LI", or even "L.LI" if no confusion arises, 
if for any bounded open set $E \subset {\mathbb R}^{n}$ and 
any $\phi \in G(E)$ we have the inequality: 
\begin{equation}
\int_{E} W(x_{0}, y_{0}, \eF \nabla \phi(y) ) \mbox{ d}y \ 
\geq \ \mid E \mid \  
W(x_{0}, y_{0}, \eF)  \ \ . 
\label{si}
\end{equation}

$W$ (or $I$) is $G$ right lower invariant ( $G$ R.LI) in  
$(x_{0}, y_{0}, \eF)$ if for any bounded open set 
$E \subset {\mathbb R}^{n}$ and any $\phi \in G(E)$ we have: 
\begin{equation}
\int_{E} W(x_{0}, y_{0},  \nabla \phi(y) \eF ) \mbox{ d}y \ 
\geq \ \mid E \mid \  
W(x_{0}, y_{0}, \eF)  \ \ .
\label{rsi}
\end{equation}
If in the relations (\ref{si}), (\ref{rsi}) we change "$\geq$"
 by "$\leq$" then we obtain the definitions 
of $G$ left upper  invariance ($G$ L.UI), respectively $G$ 
 right upper invariance ($G$ R.UI). If $W$ is right and left 
LI then we call it $G$ LI; also 
if $W$ is right (or left) lower and upper
invariant we call it right (or left) invariant. 
\label{defsi}
\end{defi}

A key observation consists in the following proposition, which
shows that quasiconvexity is a particular case of lower 
invariance.

\begin{prop}
Let us consider $(x_{0}, y_{0}, \eF) \in {\mathbb R}^{n} 
\times {\mathbb R}^{n} \times GL_{n}(R)$. Then  
$W$ is $Diff^{\infty}_{0}$ L.LI  in $(x_{0}, y_{0}, \eF)$  if
 and only if it is quasiconvex in 
the same triplet.
\label{p4}
\end{prop}

\begin{proof}
Let  $E \subset {\mathbb R}^{n}$ be an open bounded
 set and $\phi \in Diff^{\infty}_{0}(E)$. The vector field 
$\eta \ = \ \eF (\phi - id)$ verifies i) and ii) from
 definition \ref{defqc}. Therefore, if 
$W$ is quasiconvex in  $(x_{0}, y_{0}, \eF)$, we derive from 
 (\ref{qconv}) the inequality: 
$$\int_{E}W(x_{0}, y_{0}, \eF \nabla \phi(y)) \mbox{ d}
 y \ \geq \ \mid E \mid 
W(x_{0}, y_{0}, \eF) \ \ . $$ 
We implicitly used the chain of equalities $\eF + \nabla 
\eta(y) \ = \ \eF + \eF \nabla \phi (y) \ - \  \eF \ = \ 
\eF \nabla 
\phi (y)$ . We have proved that quasiconvexity implies 
$Diff^{\infty}_{0}$ L.LI

In order to prove the inverse implication we shall suppose
 that $E$ is also simply connected. 
This supposition is not restrictive according to corollary 
 3.1.1 from Ball \cite{1} (see also the references therein 
and the twin
result  contained in proposition \ref{ptwin} from this paper).  
Let us  consider  $\eta$ which satisfies  i) and ii) from 
definition \ref{defqc}. From the hypothesis upon $E$ 
the function 
$$\psi(x) = \left\{ \begin{array}{ll}
\eF x \ + \ \eta(x) & \mbox{ if } x \in E \\
\eF x & \mbox{ otherwise } \ \ , 
\end{array} \right. $$
is $C^{\infty}$ and invertible on  ${\mathbb R}^{n}$.
 We have therefore  
$\phi  \ = \ \eF^{-1} . \psi \in Diff^{\infty}_{0}(E)$ and 
$\eF \nabla \phi \ = \ \eF + \nabla \eta $. 
If $W$ is $Diff^{\infty}_{0}$ L.LI in $(x_{0}, y_{0}, \eF)$
 then  we use (\ref{si})  
with the previously defined $\phi$ in order to obtain
 (\ref{qconv}). 
\end{proof}

\begin{rk}
In  definition \ref{defsi} $G(E)$ can be replaced by
$G^{1,\infty}(E)$. This follows from the definition 
of $G^{1,\infty}(E)$ and the continuity of $W$. 
\label{rext}
\end{rk}

\begin{prop}
If $G$ is a group and  $G^{c} \ = \ G$ then  $W$ is 
$G$ L.LI in 
$(x_{0}, y_{0}, \eF) 
\subset {\mathbb R}^{n} \times {\mathbb R}^{n} \times J(G)$ 
if and only if  $W$ is $G$ R.LI in the same triplet.   
\label{bi}
\end{prop}

\begin{proof}
Let us suppose that  $W$ is  
$G$ R.LI in  $(x_{0}, y_{0}, \eF)$. We make the change
 of variable $x = \eF^{-1} y$ and 
we rewrite the hypothesis in the following way: for any 
open bounded set  
$E \subset {\mathbb R}^{n}$ and any  $\phi \in G(E)$ we 
have 
$$\int_{\eF^{-1}(E)} W(x_{0}, y_{0}, \nabla (\phi . \eF ) 
(x)) \mbox{ d} x  \ \geq \  
\mid \eF^{-1} (E) \mid W(x_{0}, y_{0}, \eF) \ \ . $$
The hypothesis of the proposition implies that the 
application  $\phi \in 
G(E) \mapsto \eF^{-1} . \phi . \eF 
\in G(\eF^{-1}(E))$ is  well defined and
 bijective. 
Therefore $W$ is  $G$ R.LI  in $(x_{0}, y_{0}, \eF)$ 
if and only if for any 
bounded open set $E \subset {\mathbb R}^{n}$ and 
for any  $\psi \in G(\eF^{-1}(E))$ we have  
$$ \int_{\eF^{-1}(E)} W(x_{0}, y_{0}, \eF \nabla \psi (x))
 \mbox{ d} x  \ \geq \  
\mid \eF^{-1} (E) \mid W(x_{0}, y_{0}, \eF) \ \ . $$
The last statement is equivalent to the fact that $W$ is  
  $G$ L.LI  in $(x_{0}, y_{0}, \eF)$. 
\end{proof}
The proposition remains true if we change
 lower invariance with upper invariance.

The following theorem shows that $G$ lower invariance of $W$ 
is a necessary condition for the existence of a minimum of 
$I(\cdot; \Omega)$ over 
$C^{1}(\Omega,{\mathbb R}^{n}) \cap W^{1,\infty}(G,
\Omega)$. 

\begin{thm} Let us suppose that there exists 
$\eu \in C^{1}(\Omega,{\mathbb R}^{n}) \cap W^{1,\infty}(G,
\Omega)$ such that
 for any  
$\phi \in G(\Omega)$, $\| \phi - id \|_{C(\Omega)} < \epsilon$,
  we have:
$$I(\eu . \phi; \Omega) \ \geq \ I(\eu; \Omega) \ \ .$$ 
$W$ is then  $G$ 
L.LI in 
$(x_{0}, \eu(x_{0}), \nabla \eu (x_{0}))$ for  any  
$x_{0} \in \Omega$ . 
\label{t1}
\end{thm}

\begin{proof} 
Let $x_{0} \in \Omega$ and  $E \subset {\mathbb R}^{n}$ 
be an open bounded set. 
We can find then  $\epsilon_{0} > 0$ such that the
 followings are true: 
\begin{enumerate}
\item $x_{0} + \epsilon_{0} \ E \ \subset \subset \Omega$ ;
\item for any  $\psi \in G(\Omega)$, $\| \psi - id \|_{C(\Omega)} 
< \epsilon_{0}$, we have:
\begin{equation}
I(\eu . \psi; \Omega) \ \geq \ I(\eu; \Omega) \ \ . 
\label{maxi}
\end{equation}
\end{enumerate}
Let us consider $\phi \in G(E)$, $\epsilon \leq \epsilon_{0}$,
 $\ef^{\epsilon}(x) \ = \ x_{0} + \epsilon x$
and $\phi^{\epsilon} \ = \ A(\ef^{\epsilon}, \phi)$. From 
lemma \ref{p1} it follows that $\phi^{\epsilon}  
\in G(\Omega)$. We can choose a sufficiently small
 $\epsilon$ such that  
$\| \phi^{\epsilon} - id \|_{C(\Omega)} < \epsilon_{0}$. 
We apply  (\ref{maxi}) with  
$\phi^{\epsilon}$ and we obtain the inequality:
$$\int_{\ef^{\epsilon}(E)} W(x, \eu(\phi^{\epsilon}(x)),
 \nabla \eu(\phi^{\epsilon}(x)) 
\nabla \phi^{\epsilon}(x)) \mbox{ d}x \ \geq \  I(\eu; 
\ef^{\epsilon}(E)) \ \ .$$
After the change of variable $\ef^{\epsilon}(y) \ = \ x$ 
the inequality becomes: 
$$\int_{E} W(\ef^{\epsilon}(y), \eu(\ef^{\epsilon}.
 \phi (y)), \nabla \eu (\ef^{\epsilon}. \phi
(y)) \nabla \phi (y) ) \ \epsilon^{n} \mbox{ d} y \ \geq \ $$ 
$$\geq \ 
\int_{E} W(\ef^{\epsilon}(y), \eu(\ef^{\epsilon}(y)), 
\nabla \eu(\ef^{\epsilon}(y))) \epsilon^{n} \mbox{ d}y \ \ 
. $$
We reduce  $\epsilon^{n}$  from the both members of the
 inequality.   
 The continuity of  $W$ and regularity of  $\eu$
 imply that when  $\epsilon$ converges to $0$ we have 
 the inequality: 
$$\int_{E} W(x_{0}, \eu(x_{0}), \nabla \eu (x_{0}) \nabla 
\phi (y)) \mbox{ d}y \ \geq \ 
\mid E \mid W (x_{0}, \eu(x_{0}), \nabla \eu(x_{0})) \ \ . $$
\end{proof}

The following proposition shows that in the definition  
\ref{defsi} the text "for any bounded open set $E$ and 
any $\phi \in G(E)$ ..." can be replaced by "there is 
a bounded open set $E$ such that for any $\phi \in G(E)$ ...".

\begin{prop}
Let us consider  $x_{0}, y_{0} \in {\mathbb R}^{n}$ and  
$\eF \in J(G)$. 
If there exists $E \subset {\mathbb R}^{n}$, bounded and 
open, such that for any   
 $\phi \in G(E)$, $\| \phi - id \|_{C(E)} < \epsilon$,  
we have 
$$\int_{E} W(x_{0}, y_{0}, \eF \nabla \phi(y) )
 \mbox{ d}y \ \geq \ \mid E \mid \  
W(x_{0}, y_{0}, \eF) $$
then $W$ is $G$ L.LI in $(x_{0}, y_{0}, \eF)$. 
\label{ptwin}
\end{prop}

\begin{proof} Let us take  $\Omega \ =  \ E$ and 
$$ I(\eu ; E)  \ = \ 
\int_{E}W(x, \eu(x), \eF \nabla \eu(x)) \mbox{ d}x \ \
 .$$ We apply theorem \ref{t1} and conclude the proof. 
\end{proof}

Any quasiconvex function $W$ is both $G$ L.LI and R.LI. This 
follows from propositions 4.1, 4.2 and the simple remark that if
$G \subset G'$ then $G'$ LI implies $G$ LI. 
\begin{opp}
Find a group $G$ and a function $W$ which is $G$ R.LI but not 
$G$ L.LI  
\label{op1}
\end{opp}

\subsection{Semicontinuity and invariance}
\indent

\begin{defi}
A functional $I: X(\Omega) \rightarrow R$ is 
left  sequentially weak* lower semicontinuous ($G$ L.LSC) in 
$\eu \in X(\Omega)$ if for any sequence
 $\phi_{h} \in G^{1,\infty}(\Omega)$  convergent to $id$ 
 we have:
$$ I(\eu) \ \leq \  \liminf_{h \rightarrow \infty  }
I(\eu . \phi_{h}) \ \ . $$
The functional $I$ is right sw*lsc ($G$ R.LSC) in $\eu$
 if for any sequence $\phi_{h} \in G^{1,\infty}(\Omega)$
 convergent to $id$  we have:
$$ I(\eu) \ \leq \  \liminf_{h \rightarrow \infty  }
I( \phi_{h} . \eu) \ \ . $$
\label{deflsc}
\end{defi}

The purpose of this section is to explore the connections
 between the  $G$ lower invariance of $W$ and the 
lower semicontinuity (in the sense of definition 
\ref{deflsc}) of the functional  $I( \cdot ; \Omega)$.
Our results generalize the ones from Morrey \cite{9}, 
Meyers \cite{3}, which show
 that quasiconvexity of $W$ is equivalent to 
lower semicontinuity (in the classical sense)  of 
$I( \cdot ; \Omega)$, if the functional is defined 
over a Sobolev vector space.

\begin{thm} 
Let $\Omega \subset {\mathbb R}^{n}$ be an open
 bounded set and 
$W: {\mathbb R}^{n} \times {\mathbb R}^{n} \times 
J(G) \rightarrow R$ continuous. 
If for any  $\psi \in G$ 
and for any sequence $\phi_{h} \in G(\Omega)$ 
 convergent 
to  $id$ we have the inequality: 
\begin{equation}
\label{slim}
I(\psi; \Omega) \ \leq \ \liminf_{h \rightarrow \infty  }
I(\psi . \phi_{h} ; \Omega)  
\end{equation} 
then 
$W$ is $G$ L.LI in any triplet of the form 
$(x, \psi(x), \nabla \psi (x)) \subset \Omega \times
\Omega  
\times J(G)$, $\phi \in G(\Omega) $. If $W \ = \ W(x, \eF)$ 
then the conclusion of the theorem is: 
$W$ is $G$ L.LI in any pair of the form 
$(x, \eF) \subset \Omega  \times J(G)$. 
\label{t2}
\end{thm}

\begin{proof} 
Let us consider $x_{1} \in \Omega$ and $h > 0$. 
 $Q_{h}$ is the cube $x_{1}^{i} < x^{i} < x^{i} + 1/h$. 
We take $\phi \in G(Q_{1})$ and $k \in N$. The extension of 
$\phi$ by periodicity over ${\mathbb R}^{n}$ 
is denoted by $\tilde{\phi}$. 
We define then: 
$$\phi_{h,k} (x) \ =  \ \left\{ \begin{array}{ll} 
( hk)^{-1} \left( \tilde{\phi}(hk(x-x_{1}) + x_{1})  -x_{1} 
\right) 
+ x_{1} & \mbox{ if } x \in 
Q_{h} \\ 
x & \mbox{ otherwise} \ \ . 
\end{array}
\right. $$
From proposition \ref{p2} and A3/ we infer that 
$\phi_{h,k} \in G$. Any set   $Q_{h}$  decomposes in  
$k^{n}$ cubes which will 
be denoted by $Q_{hk,j}$, $j = 1, ..., k^{n}$, such that 
$Q_{hk,1} = Q_{hk}$.  The corner of 
 $Q_{hk,j}$ with least distance from $x_{1}$ is denoted by 
$x_{j}$.

Let us now consider  $\psi \in G(\Omega) \cap
C^{2}(\Omega, {\mathbb R}^{n})$ and 
$y_{1}  =  \psi (x_{1})$, $\eF  =  \nabla \psi (x_{1})$. 
For a sufficiently large  
$h$ we have $Q_{h} \subset \Omega$, hence 
$I(\psi . \phi_{h,k} ; \Omega)$ makes sense. 
We decompose this integral in two parts: 
\begin{equation}
I(\psi . \phi_{h,k} ; \Omega) \ = \ I(\psi . \phi_{h,k} ; 
Q_{h}) \ + \ I(\psi  ; \Omega \setminus Q_{h}) \ \ , 
\label{prt21}
\end{equation}
$$I(\psi . \phi_{h,k} ; Q_{h}) \ =   
\sum^{k^{n}}_{j=1} \int_{Q_{hk,j}} \left[ 
W(x, \psi . \phi_{h,k}(x), \nabla ( \psi . \phi_{h,k}) (x)) 
\right. $$
\begin{equation}
  \left.  - \ W(x_{j},\psi . \phi_{h,k}(x_{j}), 
\nabla  \psi ( \phi_{h,k} (x_{j})) \nabla \phi_{h,k} (x)) 
\right] \mbox{ d}x \ + 
\label{prt22}
\end{equation}
$$  + \  \sum^{k^{n}}_{j=1} \int_{Q_{hk,j}} 
W(x_{j},\psi . \phi_{h,k}(x_{j}), 
\nabla  \psi ( \phi_{h,k} (x_{j})) \nabla \phi_{h,k} (x)) 
\mbox{ d}x \ \ . $$
Notice that $\phi_{h,k}$ converges weak* to  $id$. Because  
$W$  and $\nabla \psi$ are continuous and $\phi_{h,k}$
 converges uniformly to  $id$, it follows that 
the first sum from the right-handed member of the 
equality  (\ref{prt22}) converges to zero.

By the change of variable $y = hk(x-x_{j}) + x_{1}$
 we obtain: 
\begin{equation}
 \int_{Q_{hk,j}} W(x_{j},\psi . \phi_{h,k}(x_{j}), 
\nabla  \psi ( \phi_{h,k} (x_{j})) \nabla \phi_{h,k} (x))  
\mbox{ d}x \ = 
\label{key1}
\end{equation}
$$ = \ (hk)^{-n} \ 
\int_{Q_{1}} W(x_{j}, \psi(x_{j}), \nabla \psi (x_{j})
 \nabla \phi (y) ) \mbox{ d}y \ \ . $$
We deduce from here that the second sum of the right-handed
 member  (\ref{prt22}) is a Cauchy sum. By a passage
 to the limit as $k \rightarrow \infty$  we get the equality:   
\begin{equation}
\label{prt23}
\lim_{k \rightarrow \infty} I(\psi . \phi_{h,k} ; Q_{h}) \ = \ 
\int_{Q_{h}} \int_{Q_{1}} W(x, \psi(x), \nabla \psi(x) 
\nabla \phi(y)) \mbox{ d}y \mbox{ d}x \ \ . 
\end{equation}
From  (\ref{slim}) we have: 
$$ \liminf_{k \rightarrow \infty} I(\psi . \phi_{h,k} ; \Omega) 
 = \ \liminf_{k \rightarrow \infty} 
I(\psi . \phi_{h,k} ; \Omega \setminus Q_{h})      + 
\ I(\psi  ; \Omega \setminus Q_{h}) $$   
\begin{equation}
 \geq \ I(\psi ; Q_{h}) \ + \ I(\psi  ; \Omega \setminus Q_{h}) 
\ \ , 
\end{equation} 
therefore (\ref{prt23}) implies that: 
\begin{equation} \label{prt24} 
\int_{Q_{h}} \int_{Q_{1}} W(x, \psi(x), \nabla \psi(x) \nabla 
\phi(y)) \mbox{ d}y \mbox{ d}x  \ \geq  
\end{equation}
$$\geq \ \int_{Q_{h}} W(x, \psi(x), \nabla \psi(x) ) \mbox{ d}x 
\ \ .$$ 
We multiply the relation (\ref{prt24}) with $h^{n}$ and  pass to 
the limit as  $h \rightarrow \infty$. 
The result is: 
$$\int_{Q_{1}} W(x_{1}, \psi(x_{1}), \nabla \psi (x_{1}) \nabla 
\phi (y)) \mbox{ d}y \ 
\geq \ W(x_{1}, \psi(x_{1}) , \nabla \psi (x_{1})) $$ 
which concludes the first part of the proof.

Now, if  $W \ = \ W(x,\eF)$ then let us notice that for any 
$x_{1} \in \Omega$ and  
$\eF \in J(G)$ there exists $\psi \in G(\Omega)$ such that 
$\nabla \psi(x_{1}) \ = \ \eF$, therefore 
we can apply what we have already proved in order to obtain
 the second conclusion of the theorem. 
\end{proof}

\begin{thm} 
Let $W: {\mathbb R}^{n} \times {\mathbb R}^{n} \times J(G)
 \rightarrow R$ be a continuous function.
If $W$ is  $G$ R.LI in any triplet of the form  
$(x, y, \eF) \subset \Omega \times \Omega  
\times J(G)$,  then for  any  $\psi \in
W^{1,\infty}(G,\Omega)$ and any sequence convergent to  $id$
$\phi_{h} \in G(\Omega')$, where $\psi(\Omega) \subset \Omega'$, 
we have the inequality: 
\begin{equation}
\label{dlim}
I(\psi; \Omega) \ \leq \ \liminf_{h \rightarrow \infty} 
I( \phi_{h} . \psi ; \Omega) \ \ .  
\end{equation} 
\label{t3}
\end{thm}

\begin{proof}
For the proof  is not restrictive to consider that 
$\psi(\Omega) \subset \Omega$. 
Let  $G_{\nu}$ be the cubic lattice   constructed from the
 cube  $0 \leq x^{i} \leq 2^{-\nu}$ and let  
$\Gamma_{\nu}$ be the reunion of all cubes of  $G_{\nu}$
 included in $\Omega$. Let us consider 
 $\psi \in W^{1, \infty}
(G,\Omega)$ and a sequence $\phi_{k} \in G(\Omega)$ 
convergent to $id$. 

Fix  $\epsilon > 0$; there exists  $\nu'$, sufficiently 
large such that:
\begin{equation}
\mid I( \phi_{k} . \psi ; \Omega \setminus \Gamma_{\nu'})
 \mid \ < \ \epsilon \  \  \forall \ k \in
N  \  \  , 
\label{prt31}
\end{equation}
\begin{equation}
\mid I(\psi  ; \Omega \setminus \Gamma_{\nu'}) \mid \ < \
 \epsilon \ \ . 
\label{prt32}
\end{equation}
For any $\nu > \nu'$, with the notations from the proof of
 the theorem \ref{t2}, 
we write  $\Gamma_{\nu'}$ in the following way:  
$\Gamma_{\nu'} \ = \ \cup^{N_{\nu'}}_{h = 1} Q_{h}$.  
The integral $I( \phi_{k}. \psi; \Gamma_{\nu'})$ can be 
regarded as a sum of two terms:
$$I( \phi_{k} . \psi ; \Gamma_{\nu'}) \ = \  \int_{\Gamma_{\nu'}}
 \left[ W( x,  \phi_{k} . \psi (x) , 
\nabla ( \phi_{k}. \psi) (x) ) \ - \right.$$
\begin{equation} 
\left. - \ W ( x, \psi(x), \nabla ( \phi_{k} . \psi) (x)) 
\right] 
\mbox{ d}x 
\label{prt33}
\end{equation}
$$ + \   \int_{\Gamma_{\nu'}}  W ( x, \psi(x), \nabla
 ( \phi_{k}. \psi ) (x)) \mbox{ d}x \ \ . $$
The first integral  of the right-handed member of 
(\ref{prt33}) converges to zero as  $k \rightarrow \infty$. 

We can take in any cube $Q_{h} \in \Gamma_{\nu}$ a point 
$x_{h,\nu}$ such that it is a Lebesgue point for all $\nabla
\phi_{k}$.  
For any $\ev \in L^{1}_{loc}({\mathbb R}^{n}, {\mathbb R}^{m})$ such 
that all $x_{h,\nu}$ are Lebesgue points 
 we make the notation: 
$$ \overline{\ev}(x) \ = \ \left\{ \begin{array}{ll}
\ev(x_{h,\nu}) & \mbox{ if } x \in Q_{h} \\
x & \mbox{ otherwise } \ \ \ . 
\end{array} \right. $$
The second integral from the right-handed member of 
(\ref{prt33}) can be written as a sum 
$J_{1} + J_{2} + J_{3}$, with:
\begin{equation}
\label{prt34}
J_{1} \ = \ \int_{\Gamma_{\nu'}} \left[ W(x, \psi(x), 
\nabla  \phi_{k} (\psi(x)) \nabla  \psi  
(x)) \ - \ \right.
\end{equation}
$$\left. - \ W(\overline{x}, \overline{\psi}(x),  \nabla
 \phi_{k} (\psi(x)) 
 \overline{\nabla  \psi } (x))     \right] \mbox{ d}x \ \ , $$
\begin{equation}
\label{prt35}
J_{2} \ = \ \int_{\Gamma_{\nu'}} \left[  W(\overline{x}, 
\overline{\psi}(x),  \nabla \phi_{k} (\psi(x)) 
 \overline{\nabla  \psi } (x)) \ - \right.
\end{equation}
$$\left. - \  W(\overline{x}, \overline{\psi}(x),   
\overline{\nabla  \psi } (x)) \right] \mbox{ d}x \ \ , $$
\begin{equation}
\label{prt36}
J_{3} \ = \ \int_{\Gamma_{\nu'}}  W(\overline{x}, \overline{\psi}(x),   
\overline{\nabla  \psi } (x))  \mbox{ d}x \ \ . 
\end{equation}
From the continuity of $W$, the boundedness of 
$\nabla \psi$ and the uniform boundedness of $\nabla \phi_{k}$ 
we deduce that: 
\begin{enumerate}
\item[1.] $J_{1}$ converges to zero  uniformly with
 respect to $k$, 
\item[2.] $J_{3}$ converges to $I(\psi; \Gamma_{\nu'})$ 
\end{enumerate}
as $\nu \rightarrow \infty$.

If  $W$ is $G$ R.LI  in any triplet $(x, y \eF) 
\subset \Omega  \times \Omega \times J(G)$ then:
\begin{equation}
\liminf_{k \rightarrow \infty} J_{2} \ \geq \ 0 \ \ . 
\label{prt37}
\end{equation}
From the convergences mentioned at 1., 2.  and from the 
relations (\ref{prt33}), (\ref{prt37}) we obtain:
$$\liminf_{k \rightarrow \infty} I(\phi_{k} . \psi ; 
\Gamma_{\nu'}) \ \geq \ I(\psi ; \Gamma_{\nu'}) \ \ . $$
The latter relation, together with (\ref{prt31}), (\ref{prt32}),
 lead us to the inequality:
$$\liminf_{k \rightarrow \infty} I(\phi_{k} . \psi ; \Omega) \ 
\geq \  I(\psi; \Omega) \ - \ 2 \epsilon \ \ , \mbox{ q.e.d. } $$
\end{proof}

\begin{rk}
If  $G$ is a group then $I(\cdot ; \Omega)$ is $G$ L.LSC 
over $G$ if and only if it is $G$ R.LSC. Moreover, the 
left (or right) lower semicontinuity  over $G$ are equivalent
 with classical lower semicontinuity.   
\label{glscd}
\end{rk}

The theorem  \ref{t2} essentially says that if 
 $I(\cdot; \Omega)$ is $G$ L.LSC  then   $W$ is $G$ L.LI 
The theorem  \ref{t3} asserts that if 
  $W$ is $G$  R.LI   then  
  $I(\cdot; \Omega)$ is $G$ R.LSC. 

\begin{opp}
Are the inverse  implications true?
\label{op2}
\end{opp}
 Our guess is that they are {\it not} generally true. 
Notice that in the proof of the theorem  \ref{t2} there
 is a key equality (\ref{key1}); 
all the proof but this equality could be rewritten with
 the hypothesis that   $I$ is $G$ R.LSC and 
the conclusion  would be that  $W$ is $G$ R.LI 
Analogous remarks can be made for the theorem \ref{t3}. 
The key step in the proof of this theorem is the uniform
convergence to 
zero of the term  $J_{1}$, defined in (\ref{prt34}). 

With the results from this section the formulation of 
the open problem \ref{op1} becomes clear. Indeed, we are interested 
to find a function $W: J(G) \rightarrow R$ such that:
\begin{enumerate}
\item[i)] the integral
functional which is generated by $W$ is (left or right) lower
semicontinuous;
\item[ii)] there is no quasiconvex function $W^{*}$ with the property:
$$\forall \eu \in G^{1,\infty}(\Omega) \ \ , \ 
\int_{\Omega} W(\nabla \eu(x)) \mbox{ d}x \ = \ 
\int_{\Omega} W^{*}(\nabla \eu(x)) \mbox{ d}x \ \ . $$
\end{enumerate}
Suppose that we have found a group $G$ and a function $W$ which bare 
the open problem \ref{op1}. Then, according to proposition \ref{bi}, 
$G \not = G^{c}$. If item ii) were true, then the integral functional 
$I$ generated by $W$ would be lower semicontinuous, therefore by 
theorem \ref{t2} $W$ would be $G$ L.LI, which contradicts the
hypothesis.

\begin{opp}
Suppose that $G = G^{c}$. By remark \ref{glscd} and proposition
\ref{bi} if the  integral
functional $I$ is $G$ L.LSC  then it is $G$ R.LSC. Find 
a potential $W$ and a group $G = G^{c}$ such that the functional 
$I$ generated by $W$ is $G$ R.LSC but not $G$ L.LSC. 
\label{op3}
\end{opp}

If there is a function $W$ which responds to the open problem 
\ref{op3} then the open problem \ref{op2} would have a negative
answer.

\vspace{1.5cm}

Ball introduced in \cite{1}, definition 3.2 and theorem 3.3, the 
notion of rank one convexity.  In order to generalize this notion we 
introduce first the following definition.  

\begin{defi}
$TG(\Omega)$ is the class of all vector fields $\eta \in 
C^{\infty}_{0}(\Omega, {\mathbb R}^{n})$ such that 
the one parameter flow $\phi_{t}$, $t \in (a,b)$ with 
$a< 0 < b$, defined by 
$$\dot{\phi}_{t}(x) \ =  \ \eta(\phi_{t}(x))  \ \ \ , \ \
\phi_{0} \ = \ id $$
lies in $G$, that is 
$$\phi_{t} \in G \ \ \  \forall \ t \in (a,b) \ \ \ . $$
\label{dtg}
\end{defi}

\begin{thm}
Let $W \ = \ W(\eF)$ be  a $G$ L.LI, defined on an open 
neighbourhood of $G$. Then for any 
$\eF \in J(G)$ and for any $\eta \in TG$ we have the inequality: 
\begin{equation}
\frac{\partial^{2} W}{\partial \eF_{ij} \partial \eF_{mr}} 
(F) \eF_{ik} \eF_{mp} \ \int_{\Omega} \eta_{k,j} \eta_{p,r} 
\mbox{ d}x \ \geq \ 0 \ \ \ .
\label{legh}
\end{equation}
\label{tlegh}
\end{thm}

\begin{proof} 
 Consider the function 
$$I(t) \ =  \ \int_{\Omega} W(\eF \nabla \phi_{t}(x)) 
\mbox{ d}x \ \ . $$
This is a $C^{2}$ function which has a minimum at $t=0$, 
according to hypothesis upon $W$  (\ref{si}).  This fact implies 
that 
$$\frac{\partial I}{\partial t}(0) \ = \ 0  \ \ \ , \ \ 
\frac{\partial^{2} I}{\partial t^{2}} (0) \ \geq \ 0 \ \ . $$
The first variation of $I$ has  the form: 
$$\frac{\partial I}{\partial t}(t) \ = \ 
\int_{\Omega} \frac{\partial W}{\partial \eF_{ij}} 
(\eF \nabla \phi_{t}) 
\eF_{ik} \dot{\phi_{t}}_{k,j} \mbox{ d}x $$ 
hence for $t=0$ we obtain 
$$\frac{\partial I}{\partial t}(0) \ = \ 
\frac{\partial W}{\partial \eF_{ij}}(\eF) \eF_{ik} \int_{\Omega}
\eta_{k,j} \mbox{ d}x \ \ . $$
The integral from the right-handed member is obviously null
because $\eta$ has compact support in $\Omega$, therefore we 
obtain a trivial identity. 

The second variation of $I$ is a sum of three terms:
$$\frac{\partial^{2} I}{\partial t^{2}} (0) \ = \ A \ + \ B \ + 
\ C \ \ , $$
\begin{equation}
A \ = \ \int_{\Omega} \frac{\partial W}{\partial \eF_{ij}} (\eF) 
\eF_{ik} \eta_{k,l} \eta_{l,j} \mbox{ d} x \ \ , 
\end{equation}
\begin{equation}
B \ = \ \int_{\Omega} \frac{\partial W}{\partial \eF_{ij}} (\eF) 
\eF_{ik} \eta_{k,jl} \eta_{l} \mbox{ d} x \ \ , 
\end{equation}
\begin{equation}
C \ = \ \frac{\partial^{2} W}{\partial \eF_{ij} \partial 
\eF_{mr}} (F) \eF_{ik} \eF_{mp} \ \int_{\Omega} \eta_{k,j} 
\eta_{p,r} \mbox{ d}x  \ \ . 
\end{equation}
An integration by parts argument shows that $A + B = 0$
therefore we obtain $C \geq 0$, q.e.d. 
\end{proof}

The generalized rank-one convexity is defined further.

\begin{defi}
A $C^{2}$ function $W \ = \ W(\eF)$ is $G$ (left) 
rank one convex at $\eF \in J(G)$ if for any $\eta \in TG$
 the relation (\ref{legh}) is true. 
\label{grank}
\end{defi}

\begin{rk}
If we take $G \ = \ Diff^{\infty}_{0}({\mathbb R}^{n})$ then the 
relation (\ref{legh}) becomes 
 the Hadamard-Legendre inequality (see Hadamard  \cite{11}, Ball 
\cite{1'} and the references therein). 
Indeed, for this 
group we have $TG(\Omega) \ = \ C^{\infty}_{0}(\Omega, 
{\mathbb R}^{n})$, hence for any $\eta \in TG(\Omega)$ and 
$\eF \in J(G) = GL_{n}(R)$, the vector field $\eF \eta$ 
belongs to $TG(\Omega)$. Therefore the relation (\ref{legh}) 
can be written as: 
\begin{equation} 
\frac{\partial^{2} W}{\partial \eF_{ij} \partial 
\eF_{kl}} (F)\ \int_{\Omega} \eta_{i,j} 
\eta_{k,l} \mbox{ d}x \ \geq \ 0 
\label{leg}
\end{equation}
for any $\eta \in C^{\infty}_{0}$.  An argument from Ball 
\cite{1'}, proof of Theorem 3.4, allow us to consider piecewise 
affine vector fields $\eta$. It can be shown that  (\ref{leg}) 
implies the Legendre-Hadamard inequality: 
\begin{equation}
\frac{\partial^{2} W}{\partial \eF_{ij} \partial 
\eF_{kl}} (F)\ a_{i} a_{k} b_{j} b_{l} \ \geq \ 0  \ \ \ , 
\label{trueleg}
\end{equation}
for any vectors $a,b \in {\mathbb R}^{n}$ (see also remark 6.2). 
\label{rrank}
\end{rk}

\begin{rk}
Same arguments as in the previous remark, but for the group of 
volume-preserving diffeomorphisms $Diff^{\infty}_{0}(dx)$ show 
that the relation (\ref{legh}) implies  the following 
inequality: 
\begin{equation}
\frac{\partial^{2} W}{\partial \eF_{ij} \partial 
\eF_{mr}} (F) \eF_{ik} \eF_{mp} \ a_{k} a_{p} b_{j} b_{r} \ 
\geq \ 0  \ \ \ ,  
\label{dxleg}
\end{equation}
for any orthogonal vectors $a,b \in {\mathbb R}^{n}$, $a \cdot b \ = \ 0$. 
\label{dxrank}
\end{rk}
Theorem \ref{tlegh} has a correspondent for the case of $G$
R.LI functions. We leave this theorem to the reader.

As a corollary we have:

\begin{prop}
Any $G$ L.LI function 
$W \ = \ W(\eF)$ is $G$ rank one convex. 
\label{corrk}
\end{prop}

\begin{proof}
The result is obtained from definition \ref{grank} and theorem 
\ref{tlegh}.  
\end{proof}

%
%           LAGRANGIENI NULI
%
%
%creat pe 29.01.98, destinat includerii in full_lecture.tex
%
%
%
%
\section{Null lagrangians and group invariants}

Let us consider the pair $(X,G)$ such that $X = X . G$ ($X = G . X$ 
respectively) and a   
functional $I$  representable in integral form:
$$I(\eu) \ = \  \int_{\Omega} L(x, \eu(x), \nabla \eu (x))
\mbox{ d}x \ \ $$ 
for any $\eu \in X$. 

\begin{defi}
 A $C^{2}$function 
$L: R^{n} \times R^{n} \times J(X) \rightarrow R$ is a
$(X,G)$ invariant at left (abbreviated i.l.) if   for any $\eu \in X$, 
 $\phi \in G(\Omega)$ 
we have $I(\eu .\phi) \ = \  I(\eu)$. Right invariants are defined 
in a similar way. If $L \ = \ L(\eF)$   then we call it a 
 homogeneous $(X,G)$ left  invariant.  
\label{dinli}
\end{defi}

\begin{defi}
A $C^{2}$function 
$L: R^{n} \times R^{n} \times J(G) \rightarrow R$ is a
$G$ null lagrangian  if  for any  $\eF \in J(G)$ and  
$\phi \in G(\Omega)$ we have 
$I(\eF \phi) \ = \  I(\eF)$. 
If $L$ depends only on it's third variable then is called a
homogeneous $G$ null lagrangian. 
\label{nldef}
\end{defi}

\begin{rk}
Suppose that $G \subset X$ (equivalently $id \in X$) and $J(G)
\subset X$. 
Then any $G$ invariant is a $G$ null lagrangian. 
\label{erk1}
\end{rk}

\begin{rk}
If we take $X = J(G).G$ then $L$ is a homogeneous $(X,G)$ 
invariant at left if and only if it is a $G$ null lagrangian at left.  
\label{erk2}
\end{rk}

\begin{prop}
If $L \ = \ L(x , y , \eF)$ is a $(W^{1,\infty}(G, \Omega), G)$ 
i.l. then for any $x,y \in R^{n}$ the mapping $\eF \in J(G) 
\mapsto L(x,y, \eF)$ is a $G$ n.l.l.
\label{p3v}
\end{prop}

\begin{proof}
Direct consequence of theorem \ref{t2}. Indeed, from left continuity 
follows that $L(x,y, \cdot)$ and $-L(x,y,\cdot)$ are both $G$LLI 
\end{proof}

\subsection{Examples}

Let us consider the case $X= W^{1,\infty}(GL_{n}(R), R^{n})$ and 
$G = Diff_{0}^{\infty}(R^{n})$. By proposition \ref{p4} 
any null lagrangian at left is a classical null lagrangian. 
The class of null lagrangians is known (see Ball, Currie \& Olver \cite{21} or 
Olver \& Sivaloganathan \cite{22}). For $n=3$ for example, 
any homogeneous null lagrangian is a linear combination of 
$\eF_{ij}$, $adj \ \eF_{ij}$ and $\det \eF$. 

In the particular case that we have chosen the homogeneous null
lagrangians are also invariants at left. Indeed, take $\Omega$ simply 
connected and smooth, 
$\eu \in X$ and compute $I(\eu; \Omega)$, where $I$ is generated by a
null lagrangian. We obtain: 
$$\int_{\Omega} \det \ \nabla \eu \mbox{ d}x \ = \ \mid \eu( \Omega) 
\mid \ \ , $$
$$\int_{\Omega}  (\nabla \eu)_{ij} \mbox{ d}x \ = \ \int_{\paromega}
\eu_{i} \en_{j} \mbox{ d}s \ \ , $$
$$\int_{\Omega} adj \ (\nabla \eu)_{ij} \mbox{ d}x \ = \ 
\int_{\paromega} (\eu \wedge \en)_{ij} \mbox{ d}s \ \ . $$ 
We see that generally $I(\eu; \Omega)$ depends only on $\Omega$ and 
the value of $\eu$ on $\paromega$. These rest the same under a
composition at left with any $\phi \in Diff_{0}^{\infty}(\Omega)$.  
These considerations prove proposition \ref{pinvlag} from the
introduction. 

Let us choose now $X \ = \ Diff^{\infty}(R^{n})$ and 
$G = [SL_{n}(R)]$. The (integral) invariants given by the classical 
null lagrangians are trivial in this case. An easy $[SL_{n}(R)]$ 
invariant at left turns to be $L(\eF) \ = \ \log \det \eF$. 
Indeed, consider $\eu \in X$ and $\phi \in G(\Omega)$. We have then 
$$I(\eu . \phi^{-1};\Omega) \ =  \ \int_{\Omega} 
\left( \log \det \nabla \eu(x) \ - \ \log \det \nabla \phi(x) \right) 
\det \nabla \phi(x) \mbox{ d}x \ = \ I(\eu; \Omega) \ \ . $$
Unfortunately the restriction of $L$ to $SL_{n}(R)$ equals $0$, the 
most trivial null lagrangian.

The classical null lagrangians are all derived from the determinant. 
The determinant function $\ef \mapsto \det \eF$ transforms matrix 
multiplication into number multiplication and basically this is the
reason which makes determinant to be a null lagrangian. 

The following theorem shows an interesting change of behaviour in the
case of $G$ null lagrangians with $J(G) \subset SL_{n}$. In this
situation we find a class of non-classical null lagrangians which 
transform matrix multiplication into number addition. 

Unfortunately we have not been able to find a group $M \subset SL_{n}$ 
which has non-trivial characters and $[M]$ acts transitively on
$R^{n}$. However, the following theorem might offer an 
illustration of a general phenomenon concerning non-classical null
lagrangians.

\begin{thm}
Let $M \subset  SL_{n}(R)$ be a Lie subgroup  of the linear group 
of $n \times n$ matrices with real coefficients and positive
determinant. By a character of $M$ we mean any homeomorphism 
$\chi: M \rightarrow (0, + \infty)$ from $M$ to the multiplicative 
group of $R$. For any character $\chi$ of $M$  the function: 
\begin{equation}
W: M \rightarrow R \ , \ W(\eF) \ = \  \log \ \chi(\eF)  
\label{char}
\end{equation}
is a $[M]$ null lagrangian at left. 
\label{mainex}
\end{thm}

\begin{proof}
Let us consider $\phi \in [M](E)$, $\eF \in M$ and 
$\eta \in T[M]$, $supp \ \eta \subset \subset E$. 
We denote by $\phi_{t}$ the one parameter 
group generated by $\eta$ and we introduce the function:
$$g(t) \ = \ I(\phi_{t}; E) \ = \  \int_{E} W(\nabla \phi_{t}(x)) 
\mbox{ d}x \ \ . $$
 We have then:
$$g(t_{1} + t_{2}) 
 \ = \ \int_{E} W(\nabla \phi_{t_{2}}(\phi_{t_{1}}(x))) \mbox{ d}x \ + 
 \ \int_{E}  W( \nabla \phi_{t_{1}}(x))\mbox{ d}x \ \ . $$
By the change of variables $y = \phi_{t_{1}}(x)$ in the first integral 
 and taking account of the equality  $\det \nabla
\phi_{t} (y) = 1$ we obtain: 
$$g(t_{1}+ t_{2}) \ = \ g(t_{1}) + g(t_{2}) \ \ $$
therefore $g$ is linear. 

The function $g$ is also differentiable. We compute the first
derivative of $g$ at $t=0$ and we obtain: 
$$\frac{\partial g}{\partial t} (0) \ = \ 
\int_{E} \frac{\partial}{\partial \eF} \log \chi ( {\bf 1}) 
\nabla \eta (x) \mbox{ d} x \ = \ 0 $$ 
therefore $g$ is a constant function. From here we derive that 
$g$ is constant (provided that the exponential map covers a
neighbourhood of the identity, which is obvious). 
\end{proof}

\begin{rk} Because $[M]$ is a local group, it follows in particular 
that $[M]^{c} = [M]$. Therefore, by proposition \ref{bi}, any  
$[M]$ n.l.l. is a $[M]$ n.l.r.
\end{rk}

\subsection{Properties of null lagrangians}

\begin{thm} 
The following statements are true:
\begin{enumerate}
\item[(i)] If $I(\cdot ; \Omega)$ is left continuous then 
$W$ is a $G$ n.l.l. 
\item[(ii)] Suppose that  $W$ is a $C^{2}$ function defined 
over an open neighbourhood of $J(G)$.  
If $W$ is a $G$ n.l.l. then for any $\eF \in J(G)$ and 
 $\eta \in TG(\Omega)$ we have the equality: 
\begin{equation}
\frac{\partial^{2} W}{\partial \eF_{ij} \partial 
\eF_{mr}} (F) \eF_{ik} \eF_{mp} \ \int_{\Omega} \eta_{k,j} 
\eta_{p,r} \mbox{ d}x \ = \ 0 \ \ . 
\label{lgheq}
\end{equation}
\end{enumerate}
\label{tfir}
\end{thm}

\begin{proof}
The first statement comes from theorem \ref{t2} applied for 
$I(\cdot; \Omega)$ and $-I(\cdot; \Omega)$. The second statement 
is a consequence of the theorem \ref{tlegh}. 
\end{proof}

\begin{rk}
We particularly see that if $I(\cdot; \Omega)$ is left
invariant then $W$ is a $G$ n.l.l.  
\end{rk}

\begin{rk}
The equality (\ref{lgheq}) becomes 
the Hadamard-Legendre equality if $G =
Diff^{\infty}_{0}({\mathbb R}^{n})$.  In this case (\ref{lgheq}) is
equivalent to 
\begin{equation}
\frac{\partial^{2} W}{\partial \eF_{ij} \partial \eF_{kl}}
(\eF) \ = \ - \frac{\partial^{2} W}{\partial \eF_{il} 
\partial \eF_{kj}} (\eF) \ \ . 
\label{parhl}
\end{equation}
Generally, if for any $\eF \in J(G)$ 
we have $\eF^{T} \in J(G)$ (where $\eF^{T}$ is the transpose 
of $\eF$), then (\ref{parhl}) implies the classical Hadamard-
Legendre equality: 
$$ \frac{\partial^{2} W}{\partial \eF_{ij} \partial 
\eF_{kl}} (F) a_{i} a_{k} b_{j} b_{l} \ = \ 0 $$
for any $a, b \in {\mathbb R}^{n}$. 
\label{rpar}
\end{rk}

Let us denote by $ExpG(\Omega)$ the class of all $\phi \in 
G(\Omega)$ for which  there exist $\eta \in TG(\Omega)$ and 
$\tau \in R$ such that $\phi \ = \ \phi_{\tau, \eta}$, where 
$\phi_{s, \eta}$ is the one parameter flow generated by $\eta$.

\begin{thm}
Suppose that  $W$ is a $C^{2}$ function defined 
over an open neighbourhood of $J(G)$. Then the following 
statements are true: 
\begin{enumerate}
\item[(i)] If   $W$ satisfies (\ref{parhl}) for any $\eF \in 
J(G)$ then for any $\eu \in W^{1, \infty}(G, \Omega)$ and 
 $\eta \in TG(\Omega)$ we have the equality: 
\begin{equation}
\int_{\Omega} \frac{\partial}{\partial x_{j}} 
\left( \frac{\partial W}{\partial \eF_{ij}} 
(\nabla \eu(x))\right) \eu_{i,p}(x) \eta_{p}(x) \mbox{ d}x \ 
= \ 0 \ \ . 
\label{fvari}
\end{equation}
\item[(ii)] If $W$ satisfies (\ref{fvari}) for any 
$\eu \in W^{1, \infty}(G, \Omega)$ and  $\eta \in TG(\Omega)$ 
then we have 
\begin{equation}
I(\eu . \phi; \Omega) \ = \ I(\eu; \Omega) 
\label{expinv}
\end{equation}
for any $\eu \in W^{1, \infty}(G, \Omega)$ and 
$\phi \in ExpG(\Omega)$. In this case we say that $I(\cdot; 
\Omega)$ is exponentially $G$ invariant at left.   
\item[(iii)] If $I(\cdot; 
\Omega)$ is exponentially $G$ invariant at left then 
 for any $\eF \in J(G)$ and  $\phi \in ExpG(\Omega)$ 
we have: 
\begin{equation}
\int_{\Omega} W(\eF \nabla \phi(x)) \mbox{ d}x \ = \ \mid 
\Omega \mid \ W(\eF) \ \ . 
\label{einv}
\end{equation}
In this case we say that $W$ is exponentially $G$ 
left invariant. 
\item[(iv)] If $W$ is exponentially $G$ left invariant then 
$W$ satisfies the equality (\ref{lgheq}) for any $\eF \in 
J(G)$ and $\eta \in TG(\Omega)$. 
\end{enumerate}
\label{tlong}
\end{thm}

\begin{proof}
In order to prove (i) let us denote by $A$ the integral from
(\ref{fvari}). We have  the equality: 
$$A \ = \  \int_{\Omega} \frac{\partial^{2} W}{\partial \eF_{ij} 
\partial \eF_{kl}} (\nabla \eu(x)) \eu_{i,p}(x)\eta_{p}(x) 
\eu_{k,jl}(x) \mbox{ d}x \ \ .$$ 
From (\ref{parhl}) we see that $A = 0$. 

For (ii) let us take $\eta \in W^{1,\infty}(G, \Omega)$ and 
denote by $\phi_{t}$ the one parameter flow generated by $\eta$. 
Consider the function $g(t) \ = \ I(\eu . \phi_{t}^{-1};
\Omega)$. After some calculations based on integration by 
parts we obtain the equality: 
$$\frac{\partial g}{\partial t} (t) \ = \ 
\int_{\Omega} \frac{\partial}{\partial x_{j}} 
\left( \frac{\partial W}{\partial \eF_{ij}} 
(\nabla (\eu . \phi^{-1}_{t})(x))\right) 
(\eu . \phi^{-1}_{t}) _{i,p}(x) \eta_{p}(x) \mbox{ d}x \ \ . $$
From the hypothesis the right-handed member of the previous 
equality equals 0, therefore $g$ is a constant function and 
(\ref{expinv}) is proven. 

For (iii) take $\eu(x) = \eF x$ in (\ref{expinv}). 

For (iv) remark that the function $g$, previously 
defined, is constant. Therefore it's second variation at $t=0$ 
is null. We proceed as in the proof of 
the theorem \ref{tlegh} in order to compute this second 
variation and we finally obtain  (\ref{lgheq}). 
\end{proof}

\begin{opp}
Find the class of (homogeneous) $G$ null lagrangians. 
\label{opnul}
\end{opp}

In particular, are there any $[SL_{n}]$ non-classical null lagrangians?

\subsection{Polyconvex functions}

\begin{defi}
A function $W = W((x, y, \eF)$ is called $G$ (left or right)
 polyconvex if there is 
a continuous function $g: {\mathbb R}^{n} \times {\mathbb R}^{n} 
\times {\mathbb R}^{m} \rightarrow R$, convex in the third argument  and 
functions $w_{1}$, ..., $w_{m}$ $G$ null lagrangians (at 
left or right)  such that for any $\eF \in G$ we have 
$$W(\eF) \ = \ g(w_{1}(\eF), ... , w_{m}(\eF)) \ \ \ . $$
\label{polyc}
\end{defi}

\begin{thm}
Any $G$ (left or right) polyconvex function is $G$ LLI (or 
r.LI). Therefore, if $W$ is $G$ right polyconvex, then 
$I(\cdot; \Omega)$ is R.LSC.  
\label{tpoly}
\end{thm}

\begin{proof} For the first part of the theorem the 
proof is the same in the "right" or "left" cases. 
We apply the Jensen inequality to the function 
$g$: 
$$ \frac{1}{\mid E \mid} \ \int_{E} W(x_{0}, y_{0},
 \eF \nabla \phi(x)) \mbox{ d}x \ \geq \ 
 g(x_{0}, y_{0}, \frac{1}{\mid E \mid} \ \int_{E}
 w (\eF \nabla \phi ( x)) \mbox{ d}x ) \ \ . $$ 
The functions $w_{j}$ are null lagrangians, therefore we have: 
$$g(x_{0}, y_{0}, \frac{1}{\mid E \mid} \ \int_{E} 
w (\eF \nabla \phi ( x)) \mbox{ d}x ) \ = \ 
g(x_{0}, y_{0}, w(\eF)) \ = \ W(x_{0}, y_{0}, \eF) \ \ . $$
The second part of the theorem  is a consequence of the 
 theorem \ref{t3}. 
\end{proof}

Polyconvexity, in the classical sense, is of major interest 
because it is a local condition. Indeed, the class of 
(classical) null lagrangians was determined by 
Ericksen \cite{6} and it corresponds to the class of 
$Diff^{\infty}_{0}$ n.l.l. Therefore $W$ is polyconvex if and 
only if it has the form from definition 5.3, where 
$w_{i}(\cdot)$ are known functions (for example, if $n=3$ 
then any $w_{i}(\eF)$ is a linear combination of 
$\eF_{kl}$, $(ad \ \eF)_{kl}$ and $det \ \eF$). 

If one solves the open problem \ref{opnul},  the next problem to 
solve is the following:

\begin{opp}
Find the class of all $G$ polyconvex functions. Give sufficient 
local or global conditions for a function $W$ to be $G$ polyconvex. 
\label{opoly}
\end{opp}

Indeed, the knowledge of the class of $G$ n.l.l. functions would 
 transform the $G$ polyconvexity condition into a local 
one. In the case $G = G^{c}$ the (right) 
lower semicontinuity of $I(\cdot; \Omega)$ can be proved from 
the polyconvexity of $W$, which becomes easy to check if the 
$G$ n.l.l. functions are known.

\section{Conclusions}

We have introduced in this paper two notions: (left or right) 
lower invariance (abbreviated LI) and (left or right) 
lower semicontinuity (abbreviated LSC). These
notions describe the behaviour of integral functionals of the 
form: 
$$I(\eu;\Omega) \ = \ \int_{\Omega} W(x, \eu(x), \nabla \eu(x)) 
\mbox{ d}x $$
under inner or outer variations in a group of diffeomorphisms 
$G(\Omega)$. 

The first notion is a generalization of quasiconvexity in the 
sense of Morrey whilst the second one is weaker than the classical 
lower semicontinuity. 

For a given group $G$, the difference between $G$ LI 
and quasiconvexity of  the potential $W$ is the same as the difference 
between $G$ LSC and the classical LSC. 
For example, let us consider the case $G = [SL_{n}(R)]$, of the group 
of volume preserving diffeomorphisms. According to proposition 
\ref{bi}, in this case  $G$ left and right LI are
equivalent.  Theorem \ref{t2}  shows  that $G$ left 
LSC implies left LI hence right 
LI; by theorem \ref{t3} right LI implies 
right LSC We conclude that for general complete
groups $G$ right LSC is weaker than $G$ left LSC.

 If $G$ right LSC would be 
equivalent to classical LSC then  (at least for our
example) $G$ LI would be equivalent to quasiconvexity. 
Indeed, as earlier we use theorem \ref{t3} to deduce that 
right LI implies right LSC; by hypothesis 
right LSC implies classical LSC; 
classical LSC implies 
quasiconvexity. Therefore right LI implies
quasiconvexity. The inverse implication is always true, by proposition 
\ref{p4}. 

$G$ left or right LSC is in fact classical 
LSC when we restrict the class of admissible 
sequences to sequences obtained by inner or outer variations. 
Therefore any potential which generates integral functionals which 
are left or right  but not classical lsc gives an indications about 
the local behaviour of the group $G$.

In this paper we address several open problems. For reader's
convenience we rewrite them here: 

\begin{enumerate}
\item[{\bf OP1:}] Find a group $G$ and a function $W$ which is $G$ right LI but not 
$G$ left LI  
\item[{\bf OP2:}] Does $G$ left LI generally imply $G$ left LSC ?
Does $G$ right LSC generally imply $G$ right LI ?
\item[{\bf OP3:}] Find a potential $W$ and a group $G = G^{c}$ such that the functional 
$I$ generated by $W$ is $G$ right LSC but not $G$ left LSC 
\item[{\bf OP4:}] Find the class of (homogeneous) $G$ null lagrangians. 
In particular, are there any $[SL_{n}]$ non-classical 
null lagrangians?
\item[{\bf OP5:}] Find the class off all $G$ polyconvex functions. Give sufficient 
local or global conditions for a function $W$ to be $G$ polyconvex. 
\end{enumerate}

\end{document}